# The hyper-Wiener index of one-pentagonal carbon nanocone

M.H.Khalifeh, M.R.Darafsheh*, Hassan Jolany

School of Mathematics, College of Science, University of Tehran, Tehran .Iran.

E-mail addresses: khalife@khayam.ut.ac.ir(M.H.Khalifeh), darafsheh@ut.ac.ir(M.R.Darafsheh),
Hassan.jolany@khayam.ut.ac.ir (Hassan Jolany)

*Abstract:*

one-pentagonal carbon nanocone consists of one pentagone as its core surrounded by layers of hexagons .if there are *n* layers ,then the graph of this molecules is denoted by $G_n$ .In This paper our aim is to calculate the hyper-Wiener index of $G_n$ explicitly .

## 1. Introduction and Preliminary results

Let $G$ be a simple connected graph. The vertex set and edge set of $G$ are denoted by $V(G)$ and $E(G)$ respectively. The distance between two vertices $u$ and $v$ in $V(G)$ is denoted by $d_G(u,v)$ and is equal to the length of the shortest path from $u$ to $v$.

The Wiener index of $G$ is denoted by $W(G)$ and is defined by $W(G) = \sum_{\{u,v\}\subseteq V(G)} d_G(u,v)$. The name Wiener index is usual for chemical graphs, that is the graph of a chemical molecule, because Harold Wiener was the first person who considered this invariant for a Chemical graph [1].wiener used this index only for acyclic molecules in a slightly different way. But the definition of the Wiener index in terms of distances between vertices of a graph for the first time was given by Hosoya in [2].The Wiener index of graphs is extensively studied in [3], [4] and [5].In [6] new method are invented to calculate the Wiener index of a graph. a generalization of the concept of the Wiener index in recent chemical studies is the hyper-Wiener index.

Which was put forward in [7] and since then many research papers are written about this index. For example we can refer the reader to [8], [9] and [10]. We will define the hyper-wiener index of a graph $G$ as follow

**Definition1.** Let $G$ be a simple connected graph with vertex set $V(G)$. For a real number $\lambda$ define

$$W_\lambda(G) = \sum_{\{u,v\}\subseteq V(G)} d_G(u,v)^\lambda$$

Then $W_1(G)$ is the Wiener index of $G$, and the hyper-Wiener index of $G$ which is denoted by $WW(G)$ is defined by

$$WW(G) = \frac{1}{2}\sum_{\{u,v\}\subseteq V(G)} \left((d_G(u,v)^2 + d_G(u,v)\right) = \frac{1}{2}\left(W_2(G) + W_1(G)\right)$$

Some more concepts are needed for our further work which are defined below.

**Definition2.** Let $G$ be a simple connected graph with vertex set $V(G)$. For $F, K \subseteq V(G)$ the following quantity is defined:

$$D_G{}^\lambda(F, K) = \sum_{u \in F}\sum_{v \in K} d_G(u,v)^\lambda$$

Where $\lambda$ is a real number. If $K = V(G)$, then we set $D_G^\lambda(F, V(G)) = D^\lambda(F, G)$

**Definition 3.** Let $G$ be a simple connected graph. The subgraph $H$ of $G$ is called isometric and is written $H \ll G$ if $d_G(u,v) = d_H(u,v)$ for all $u, v \in V(H)$.

**Definition4.** Let $G$ be a simple connected graph, the subgraph $H$ of $G$ is called convex if it contains all the shortest paths in $G$ between each pair of its vertices.

Referring to the above concepts if $H \ll G$, then it is evident that $D_H^\lambda(V(H), V(H)) = 2W_\lambda(H)$ and if $\{V_k\}_{k=1}^n$ is a partition of $V(G)$, then

$$W_\lambda(G) = \frac{1}{2}\sum_{i=1}^n\sum_{j=1}^n D_G^\lambda(V_i, V_j)$$

To state the next result we need some explanation. Let $G$ be a simple connected graph. If $e \in E(G)$, then $G - e$ stands for the graph remaining from $G$ by deleting $e$ not its ends. Similarly if $F \subseteq E(G)$, then $G - F$ is defined to be the graph remaining from $G$ by deleting all the edges in $F$ not the end vertices.

**Theorem1.** Let $G$ be a simple connected graph, if $\{F_i\}_{i=1}^n$ is a partition of $E(G)$ such that each of $(G - F_i)$ is a graph with two convex connected components $G_1^i$ and $G_2^i$, then

$$W_{\lambda+1}(G) = nW_\lambda(G) - \sum_{i=1}^{n}\left(W_\lambda(G_1^i) + W_\lambda(G_2^i)\right)$$

Where $\lambda$ is a real number.

**Proof.** See [11].

**Theorem 1.** Let $G$ be a simple connected graph, if $\{F_i\}_{i=1}^{n}$ is a partition of $E(G)$ such that each of $(G - F_i)$ is a graph with two convex connected components $G_1^i$ and $G_2^i$, then the Wiener index of $G$ is:

$$W(G) = \sum_{i=1}^{n}|V(G_1^i)|\,|V(G_2^i)|$$

**Proof.** See [11]

## 2. Computing with subgraphs

As we mentioned earlier our aim is to calculate the hyper-Wiener index of the one-pentagonal carbon nanocone. The graph of this molecule consist one pentagone surrounded by layers of hexagons. If there are n layers, then this graph is denoted by $G_n$. In the following the graph of $G_6$ is drawn:

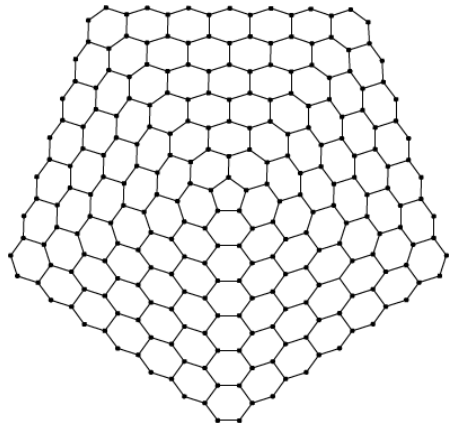

Figur1: The graph of $G_6$

Our calculations are based on Theorems 1 and 2. For the starting point we calculate the Wiener index of the auxiliary graphs which are denoted by $Z_{n,k}, M_{n,k}, Z_{n,k,l}$ and $A_n$. First we explain the above graph and draw them in special cases of $n, k$ and $l$.

$\mathbf{Z_{n,k}}$ : This graph consists of $k$ rows of hexagons with exactly $n$ hexagons in each row with two extra edges. In the following $Z_{9,5}$ is drawn:

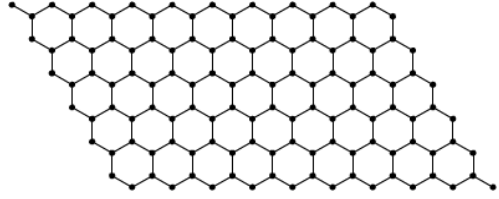

Figure 2: The graph of $Z_{9,5}$

It can be seen that the number of vertices $Z_{n,k}$ is equal to $|V(Z_{n,k})| = z_{n,k} = 2(n+1)(k+1)$.

$M_{n,k}$ : This graph consists of $k$ rows of hexagons such that in the last row ($k$ th row) there are exactly $n$ hexagons with two extra edges. In the following we draw $M_{11,6}$

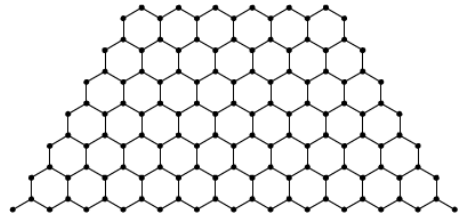

Figure 3: The graph of $M_{11,6}$

The number of vertices of this graph is equal to: $m_{n,k} = (k+1)(2n-k+3)$

$A_n$ : This graph consists of $n$ rows of hexagons such that in row $n$ there are exactly $n$ hexagnes with three extra edges. In the following we draw $A_7$

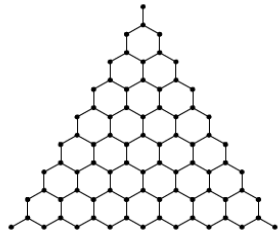

Figure 4: The graph of $A_7$

The number of vertices of the graphs is $a_n = (n+2)^2$ .

$Z_{n,k,l}$ : This graph consists of $k$ rows of hexagons such that in the first $l$ row with exactly $n$ hexagones, and from $l+1$ th to $k$ th row the number of hexagones decrease by one in each row from its previous one. We draw $Z_{11,6,2}$ in the following

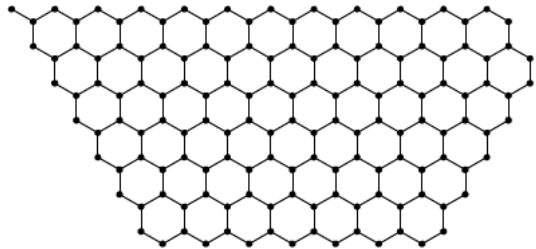

Figure 5: The graph of $Z_{11,6,2}$.

The number of vertices of this graph is: $z_{n,k,l} = 2(n+1)(k+1) - (k-l+1)^2$.

**Remark1.** For each of the graphs($Z_{n,k}$, $M_{n,k}$, $A_n$ and $Z_{n,k,l}$),in The present figure for example for the graph $Z_{n,k,l}$ ,the edges lying on each line in the following figure is equivalent to one of $F_i$ in theorems 1 and 2. Moreover they constitute a partition on the set of edges. Also for other mentioned graph similar partition exists.

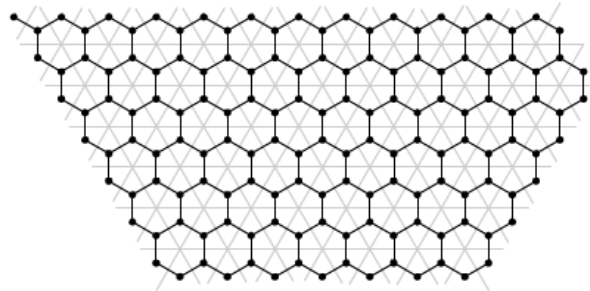

Figure 6: Partition of Edges

In figure 6 we have removed the set of edges lying on one of lines in figure 7, so by the remark 1 we have a two component graph with each component a convex subgraph of $Z_{n,k,l}$. For more explanation see [11] and [12].

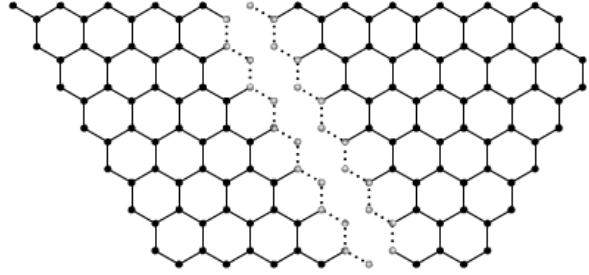

Figure 7: removing a set of edges from the graph

**Remark 2:** In theorems 1 and 2 we can see that $E(G_1^i)$, $E(G_2^i)$ and $F_i$ partition $E(G)$ and $V(G_1^i)$ and $V(G_2^i)$ partition $V(G)$ for each $i$.

**Lemma 1.** we have

$$W(A_n) = 9 + \frac{261}{10}n + 29n^2 + \frac{31}{2}n^3 + 4n^4 + \frac{2}{5}n^5$$

**Proof:** If in $A_n$ we delete a collection of edges That seem to be parallel in one row according to remark 1, then $A_n$ is partitioned into two subgraphs of the form $A_r$ and $M_{s,t}$. Therefore there is a partition of edges of the graph satisfying Theorem 2. Hence by Theorem 2 we can write

$$W(A_n) = \sum_{i=1}^{n} |V(G_1^i)| \, |V(G_2^i)| = \sum_{i=0}^{n} (a_{i-1} \times m_{n,n-i})$$

$$= 9 + \frac{261}{10}n + 29n^2 + \frac{31}{2}n^3 + 4n^4 + \frac{2}{5}n^5$$

**Lemma 2:** we have

$W(Z_{n,k}) = 1 + \frac{56}{15}k + \frac{32}{3}nk + \frac{11}{3}n + 4k^2 + \frac{1}{3}nk^4 + 6n^2k^2 + \frac{2}{3}n^2k^3 - \frac{1}{15}k^5 + \frac{4}{3}n^3 + \frac{4}{3}n^3k^2 + \frac{8}{3}n^3k + \frac{28}{3}nk^2 + \frac{8}{3}nk^3 + 4n^2 + \frac{4}{3}k^3 + \frac{28}{3}n^2k$

Proof .By referring to the graph of $Z_{n,k}$ we observe that the set of edges that seem to be parallel forms a partition of $E(Z_{n,k})$ and if we delete these edges we obtain a graph with two connected components which are convex subgraphs. Now using Theorems 2 and the fact that $V(G) = V(G_1^i) \cup V(G_2^i)$ , $1 \leq i \leq k$ , and one of the components is a graph of type $A_r$, we can write :

$$W(Z_{n,k}) = 2 \sum_{i=1}^{k} a_{i-1} \cdot (z_{n,k} - a_{i-1}) + \sum_{i=0}^{n-k-2} m_{k+i,k} \cdot (z_{n,k} - m_{k+i,k})$$

$$+\sum_{i=1}^{n} z_{n-i,k} \cdot (z_{n,k} - z_{n-i,k}) + \sum_{i=1}^{k} z_{n,k-i} \cdot (z_{n,k} - z_{n,k-i})$$

$$= 1 + \frac{56}{15}k + \frac{32}{3}nk + \frac{11}{3}n + 4k^2 + \frac{1}{3}nk^4 + 6n^2k^2 + \frac{2}{3}n^2k^3 - \frac{1}{15}k^5 +$$
$$\frac{4}{3}n^3 + \frac{4}{3}n^3k^2 + \frac{8}{3}n^3k + \frac{28}{3}nk^2 + \frac{8}{3}nk^3 + 4n^2 + \frac{4}{3}k^3 + \frac{28}{3}n^2k$$

Next we compute the Wiener index of the graph $M_{n,k}$. But this is done using the same method as above and the calculation is as follows:

$$W(M_{n,k}) = 2\sum_{i=0}^{k} a_{i-1} \cdot (m_{n,k} - a_{i-1}) + \sum_{i=0}^{n-k-1} m_{k+i,k} \cdot (m_{n,k} - m_{k+i,k})$$
$$+ \sum_{i=1}^{k} m_{n-i,k-i} \cdot (m_{n,k} - m_{n-i,k-i})$$
$$= 4 + \frac{26}{3}n + \frac{109}{15}k + \frac{1}{2}k^2 + 16nk + 4nk^2 - \frac{8}{3}k^3n - \frac{5}{2}k^3 + \frac{34}{3}n^2k + 6n^2 + \frac{2}{3}nk^4 +$$
$$\frac{8}{3}n^3k + 4n^2k^2 - \frac{4}{3}k^3n^2 + \frac{4}{3}n^3k^2 - \frac{4}{15}k^5 + \frac{4}{3}n^3$$

With the same method the Wiener index of the graph $Z_{n,k,l}$ is computed in general. the detail of the calculation is as follows:

$$W(Z_{n,k,l}) = \sum_{i=0}^{k} a_{i-1} \cdot (z_{n,k,l} - a_{i-1}) + \sum_{i=0}^{n-k-2} m_{k+i,k} \cdot (z_{n,k,l} - m_{k+i,k})$$
$$+ \sum_{i=1}^{l} m_{k-i,l-i} \cdot (z_{n,k,l} - m_{k-i,l-i}) + \sum_{i=1}^{k-l} z_{n,k-i,l} \cdot (z_{n,k,l} - z_{n,k-i,l})$$
$$+ \sum_{i=1}^{l} z_{n,l-i} \cdot (z_{n,k,l} - z_{n,l-i}) + \sum_{i=0}^{n-k+l-1} z_{k,i} \cdot (z_{n,k,l} - z_{k,i})$$
$$+ \sum_{i=0}^{k-l-1} z_{n-k+l+i,k,k-i} \cdot (z_{n,k,l} - z_{n-k+l+i,k,k-i})$$

$$= \frac{4}{5}l - \frac{1}{15}k + \frac{2}{3}n + \frac{4}{3}nk - \frac{5}{6}k^2 + \frac{13}{3}ln + \frac{4}{3}kl + \frac{3}{2}l^2 - \frac{1}{6}k^2l + \frac{23}{6}kl^2 + 2k^2l^2 - \frac{7}{6}k^3 -$$
$$\frac{2}{3}k^4 + 7kln + 4l^2nk - \frac{4}{3}l^3nk + \frac{1}{3}l^2n + \frac{10}{3}n^2k + 2n^2 - \frac{7}{6}l^3 - \frac{4}{3}l^3n - 2l^2n^2 + 4ln^2 +$$
$$8ln^2k + 4k^2l^2n - \frac{8}{3}k^3ln + 4k^2n^2l - 2n^2l^2k - \frac{4}{3}l^3k + \frac{2}{3}nk^4 + \frac{2}{3}k^4l - \frac{2}{3}k^3l^2 -$$
$$\frac{4}{3}k^3n^2 + \frac{1}{3}kl^4 - \frac{4}{15}k^5 + \frac{4}{3}n^3 - \frac{2}{15}l^5 + \frac{4}{3}n^3k^2 + \frac{8}{3}n^3k - \frac{1}{3}l^4n$$

Now by Theorem1, if the graph $G$ has the appropriate property, then for the hyper-Wiener index of $G$ we can write:

$$WW(G) = \frac{n+1}{2}W(G) - \frac{1}{2}\sum_{i=1}^{n}\left(W(G_1^i) + W(G_2^i)\right)$$

Therefore to find the hyper-Wiener index of $G$ one has to find the Wiener index of the subgraphs $G_1^i$ and $G_2^i$ instead of $|V(G_2^i)|$ and $|V(G_2^i)|$ and the Wiener index of $G$ ,but $G_1^i$ and $G_2^i$ are isomorph to the graphs $Z_{n,k}$ or $M_{n,k}$ or $A_n$ or $Z_{n,k,l}$. Therefore by our previous Computation one can see that:

$$WW(Z_{n,k}) = (n+k+1)W(Z_{n,k})$$
$$-\frac{1}{2}\Bigg[\sum_{i=0}^{k} 2\Big(W(A_{i-1}) + W(Z_{n,k,k-i})\Big)$$
$$+ \sum_{i=0}^{n-k-2} \Big(W(M_{k+i,k}) + W(M_{n-2-i,k})\Big) + 2\Big(W(A_{k-1}) + W(M_{n-1,k})\Big)$$
$$+ \sum_{i=1}^{k} \Big(W(Z_{n,k-i}) + W(Z_{n,i-1})\Big) + \sum_{i=1}^{n} \Big(W(Z_{k,n-i}) + W(Z_{k,i-1})\Big)\Bigg]$$

$$= 1 + \frac{43}{10}k + \frac{21}{5}n + \frac{211}{36}k^2 + \frac{1043}{180}n^2 + \frac{1283}{90}nk + 2nk^4 + \frac{38}{3}n^2k^2 + \frac{10}{3}k^3n^2$$
$$+ \frac{13}{3}k^2n^3 + \frac{67}{9}kn^3 + \frac{47}{3}nk^2 + \frac{67}{9}nk^3 + \frac{47}{3}kn^2 - \frac{2}{15}k^5 + \frac{10}{3}n^3 + \frac{10}{3}k^3$$
$$+ \frac{3}{2}kn^4 + \frac{25}{36}k^4 + \frac{1}{2}k^4n^2 + \frac{1}{5}k^5n + \frac{2}{9}k^3n^3 - \frac{1}{30}n^5 + \frac{25}{36}n^4 - \frac{1}{18}k^6$$
$$+ \frac{5}{6}k^2n^4 - \frac{1}{15}n^5k + \frac{1}{90}n^6$$

And similarly:

$$WW(M_{n,k})$$
$$= \frac{(2n+k+3)}{2}W(M_{n,k})$$
$$-\frac{1}{2}\Bigg[\sum_{i=0}^{k-1} 2\Big(W(A_{i-1}) + W(Z_{n-i,k,i+1})\Big)$$
$$+ \sum_{i=0}^{n-k-1} 2\Big(W(Z_{k,i}) + W(M_{n-i-1,k})\Big) + 2\Big(W(A_{k-1}) + W(Z_{n-k,k})\Big)$$
$$+ \sum_{i=1}^{k} \Big(W(M_{n-i,k-i}) + W(M_{n,i-1})\Big)\Bigg]$$

$$= 5 + \frac{113}{12}k + \frac{193}{15}n + \frac{4}{3}k^2n^4 - \frac{2}{15}kn^5 - \frac{34}{15}nk^5 + \frac{727}{30}nk + \frac{35}{6}nk^2 - \frac{23}{6}k^3n$$
$$+ 22n^2k - \frac{1}{3}nk^4 + \frac{26}{3}n^3k + \frac{41}{6}n^2k^2 - \frac{2}{3}k^3n^2 + 2n^3k^2 + \frac{7}{24}k^4 + \frac{2123}{180}n^2$$
$$+ \frac{133}{120}k^2 - \frac{29}{12}k^3 + \frac{14}{3}n^3 + \frac{5}{3}n^4k + \frac{1}{90}n^6 + \frac{3}{5}k^6 + \frac{10}{3}k^4n^2 - \frac{8}{3}k^3n^3 - \frac{1}{30}n^5$$
$$+ \frac{25}{36}n^4$$

### 3. Computation with one-pentagonal Carbon nanocone.

As we mentioned earlier, our aim is to calculate the hyper-Wiener index of the graph $G_n$ which is the graph of one-pentagonal carbon nanocone. It consists of a pentagon as its center and is surrounded by $n$ layers of hexagons. The graph of $G_n$ can be seen in section2. In this section using previous results we reach to our goal.

For a general graph $G$ and a subset $F \subseteq V(G)$ let us define $\langle F \rangle_G$ as the induced or generated subgraph by $F$ whose vertex set is $F$ and the edge set is:

$$E(\langle F \rangle_G) = \{uv = e \in E(G) | u, v \in F\}$$

Now consider the graph of $G_n$ and partition it in five sets $F_i$, $1 \leq i \leq 5$. For simplicity we show these partitions for $G_6$ as follows:

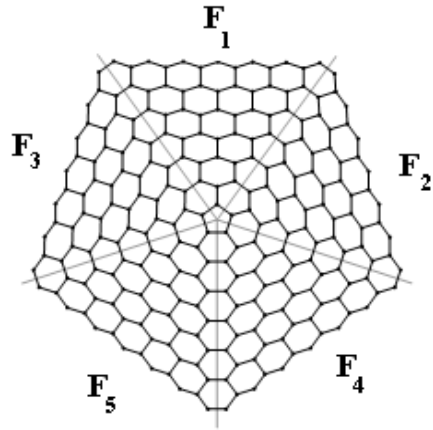

Figure 8: Partitions of $G_6$

We can see that

$$\langle F_1 \rangle_{G_n} \ll G_n$$

$$\langle F_1 \cup F_2 \rangle_{G_n} \ll G_n$$

$$\langle F_1 \cup F_2 \cup F_3 \rangle_{G_n} \ll G_n$$

We will draw the above graphs for the special case of $G_6$ :

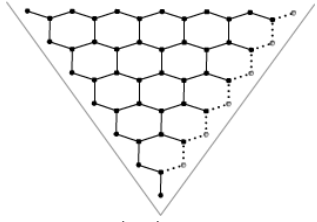 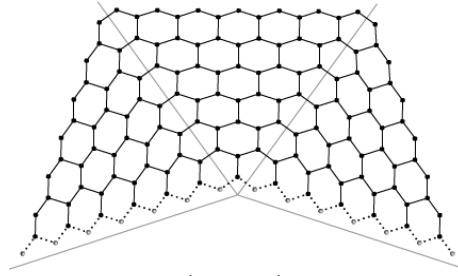 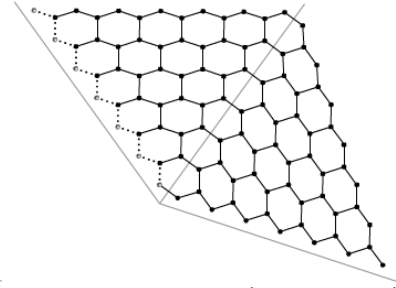

Figure 9: $\langle F_1 \rangle_{G_6}$    Figure 10: $\langle F_1 \cup F_2 \rangle_{G_6}$    Figure 11: $\langle F_1 \cup F_2 \cup F_3 \rangle_{G_6}$

Now if we consider $G_n$, we see that

$$\langle F_1 \rangle_{G_n} \cong A_n$$

$$\langle F_1 \cup F_2 \rangle_{G_n} \cong Z_{n,n}$$

$$\langle F_1 \cup F_2 \cup F_3 \rangle_{G_n} \cong M_{2n,n}$$

Theorem 3. Let $Z_{n,k}$, $M_{n,k}$ and $G_n$ be the graphs that were mentioned previously. Then for any real number λ we have :

$$W_\lambda(G_n) = 5\left(W_\lambda(M_{2n,n}) - W_\lambda(Z_{n,n})\right)$$

Proof: we saw that $\{F_i\}_{i=1}^n$ is a partition of the vertex set of $G_n$. By definition 2 we can write:

$$D^\lambda(F_1, G_n) = \sum_{i=1}^{5} D^\lambda_{G_n}(F_1, F_i) = D^\lambda_{G_n}(F_1, F_1) + 2D^\lambda_{G_n}(F_1, F_2) + 2D^\lambda_{G_n}(F_1, F_3)$$

And

$$D^\lambda_{G_n}(F_1 \cup F_2, F_1 \cup F_2) = 2D^\lambda_{G_n}(F_1, F_1) + 2D^\lambda_{G_n}(F_1, F_2)$$

And

$$D^\lambda_{G_n}(F_1 \cup F_2 \cup F_3, F_1 \cup F_2 \cup F_3) = 4D^\lambda_{G_n}(F_1, F_2) + 3D^\lambda_{G_n}(F_1, F_1) + 2D^\lambda_{G_n}(F_1, F_3)$$

Hence

$$D^\lambda(F_1, G_n) = D^\lambda_{G_n}(F_1 \cup F_2 \cup F_3, F_1 \cup F_2 \cup F_3) - D^\lambda_{G_n}(F_1 \cup F_2, F_1 \cup F_2)$$

But we already seen:

$$\langle F_1 \cup F_2 \cup F_3 \rangle_{G_n} \cong M_{2n,n} \ , \ \langle F_1 \cup F_2 \rangle_{G_n} \cong Z_{n,n}$$

And using the results

$$D^\lambda(F_1 \cup F_2 \cup F_3, F_1 \cup F_2 \cup F_3) = 2W_\lambda(M_{2n,n})$$

$$D^\lambda(F_1 \cup F_2, F_1 \cup F_2) = 2W_\lambda(Z_{n,n})$$

And finally

$$D^\lambda(F_1, G_n) = W_\lambda(M_{2n,n}) - W_\lambda(Z_{n,n})$$

Now using definition 2 we have

$$W_\lambda(G_n) = \frac{1}{2}\sum_{i=1}^{5}\sum_{j=1}^{5} D^\lambda_{G_n}(F_i, F_j) = \frac{1}{2}\sum_{i=1}^{5} D^\lambda(F_i, G_n) = \frac{5}{2} D^\lambda(F_i, G_n) \ \ 1 \le i \le 5$$

So

$$W_\lambda(G_n) = \frac{5}{2} D^\lambda(F_1, G_n)$$

And finally

$$W_\lambda(G_n) = 5\left(W_\lambda(M_{2n,n}) - W_\lambda(Z_{n,n})\right)$$

**Theorem 4.** We have

$$WW(G_n) = 20 + \frac{533}{4}n + \frac{8501}{24}n^2 + \frac{5795}{12}n^3 + \frac{8575}{24}n^4 + \frac{409}{3}n^5 + 21n^6$$

Proof. Using theorem 3 and the following fact

$$WW(G_n) = 5\left(WW(M_{2n,n}) - WW(Z_{n,n})\right)$$

Since we already calculated $WW(M_{n,k})$ and $WW(Z_{n,k})$, therefore we can calculate $WW(M_{2n,n})$ and $WW(Z_{n,n})$, as follows

$$WW(Z_{n,n}) = 1 + \frac{17}{2}n + \frac{1166}{45}n^2 + 38n^3 + \frac{521}{18}n^4 + 11n^5 + \frac{74}{45}n^6$$

$$WW(M_{2n,n}) = 5 + \frac{703}{20}n + \frac{34831}{360}n^2 + \frac{1615}{12}n^3 + \frac{7229}{72}n^4 + \frac{574}{15}n^5 + \frac{263}{45}n^6$$

Now using theorem 3 we finally obtain:

$$WW(G_n) = 20 + \frac{533}{4}n + \frac{8501}{24}n^2 + \frac{5795}{12}n^3 + \frac{8575}{24}n^4 + \frac{409}{3}n^5 + 21n^6$$

**Acknowledgement:** The authors would like to thank the research Council of the college of science through the grant no#.........